\documentclass[11pt]{article}
\usepackage{latexsym}
\usepackage{graphicx}
\usepackage{amsfonts}
\usepackage[scanall]{psfrag}
\newtheorem{theorem}{Theorem}

\newtheorem{lemma}{Lemma}
\newtheorem{proposition}{Proposition}
\newtheorem{definition}{Definition}
\newtheorem{corollary}{Corollary}
\newtheorem{notation}{Notation}
\newtheorem{remark}{Remark}

\newtheorem{hyp}{Hypothesis}
\def\C{\mathbb C}
\def\R{\mathbb R}
\def\Z{\mathbb Z}

\def\qed{$\Box\\$\medskip}
\def\phi{\varphi}
\def\epsilon{\varepsilon}

\def\proof{\medskip\noindent{\bf Proof. }}

\def\PPPF{Principal Poincar\'e Pontryagin Function }

\begin{document}

\title{\Large{\bf  \PPPF associated to some families of  Morse real polynomials}}

\author{M. Pelletier \footnote{Institut de Math\'ematiques de
Bourgogne, U.M.R. 5584 du C.N.R.S., Universit\'e de Bourgogne, B.P.
47870, 21078  Dijon Cedex, France.}, M.Uribe \footnote{Departamento
de Matem\'atica y F\'isica Aplicadas, Facultad de Ingenier\'ia,
Universidad Cat\'olica de la Ssma. Concepci\'on, Casilla 297,
Concepci\'on, Chile.}\footnote{email: michele.pelletier@normalesup.org}}

\maketitle

\begin{abstract}  It is known that the \PPPF is generically an Abelian integral. We give a sufficient condition on monodromy to ensure that it is an Abelian integral also in non generic cases. 

 In non generic cases it is an iterated integral. Uribe \cite{U,U2} gives in a special case a precise description of the Principal Poincar\'e Pontryagin Function, an iterated integral of length at most 2, involving logarithmic functions with only one ramification at a point at infinity. We extend this result to some non isodromic families of real Morse polynomials.

\end{abstract}

{\bf Keywords}: Perturbation, First return map, Iterated integrals, Monodromy, Stratification.
\\
\\
{\bf MSC}: 34M35; 34C08;14D05

\section{Introduction}
Throughout the paper $F$ denotes a Morse polynomial  $F(x,y):\C^2\rightarrow \C$ with real coefficients, of degree $d\geq 3$  and with $d$ distinct real points at infinity. It always has $(d-1)^2$ critical points but in non generic cases it can have less than $(d-1)^2$ critical values. The one-form $dF$ defines a foliation of $\C^2$. We consider a family of ovals $\delta(t)$ in regular fibers $F=t$ for $t$ in some open interval and a transverse section to these ovals, parametrized by $t$. Let $\omega(x,y)$ be a real polynomial one-form. To a one-parameter foliation defined by the perturbation  $dF+\varepsilon \omega$ for a small parameter $\varepsilon$ is associated  the displacement map which is the difference of the first return map and the identity. It is analytic with respect to $\varepsilon$. The family of ovals is destroyed if and only if the expansion with respect to $\varepsilon$ of the displacement map is not identically 0. In order to control the number of isolated zeroes of the displacement map, or in other words the number of limit cycles of the perturbed foliation, it is crucial to know the nature of the first nonzero coefficient of its expansion in $\varepsilon$. It is called the Generating Function in \cite{GI}. Following \cite{G} we will call it the Principal Poincar\'e Pontryagin Function.

It is known that it is an iterated integral \cite{Ch} and that its length depends on the monodromy group of the Milnor fibration associated to the non perturbed polynomial $F$ \cite{GI}. Generically, that is if all $(d-1)^2$ critical values are distinct, this monodromy acts transitively on the homology with complex coefficients of regular fibers $F^{-1}(t)$ and the \PPPF is an Abelian integral. It may also be an Abelian integral in non generic cases  \cite{JMP1, JMP2}. If it is not an Abelian integral, the simplest case is the one where it is a length 2 iterated integral, for example if $F$ is a triangle \cite{I}, or more generally if $F$ is the product of $d$ linear factors, $F=\ell_1.\cdots.\ell_d$, satisfying to the following Hypothesis  \cite{U, U2}.

\begin{hyp} \label{U2gen}  
The $d$ points at infinity $\ell_k=0$ are distinct, all critical points are Morse points, and the 0-level is the only critical level containing more than one critical point. The $d(d-1)/2$ intersection points of the line  $\ell_k=0$ are real.
\end{hyp}
These properties ensure that the 0-level of $F$ is what A'Campo calls a divide in  \cite{AC, AC2}, see Section 3 for the definition. We keep this terminology.

\begin{definition} A polynomial $F=\ell_1.\cdots.\ell_d$ is a generic divide in lines if it satisfies to Hypothesis \ref{U2gen}.\end{definition}
In \cite{U2} one of us proves that for generic divides in lines the \PPPF is an iterated integral of length at most 2. The proof uses  monodromy properties of  divides. The divide shows all the homology of regular fibers $F^{-1}(t)$ and allows also to compute the monodromy. Since Hypothesis \ref{U2gen} is stable it is natural to hope that the \PPPF remains a length 2 iterated integral after a small perturbation. In Section 2 we give two examples of one-parameter small perturbations of  generic divides in lines and we check that the \PPPF is still of length at most 2. Therefore we note that the fibration defined by perturbed polynomial has more monodromy operators than the fibration defined by the generic divide in lines. We also show that if the orbit of some oval  generates a codimension 1 subspace of the homology then the \PPPF is an Abelian integral.

In Section 3  we generalize examples of Section 2 by introducing

\begin{definition} A connecting family is a continous family of Morse polynomials $F_\lambda, \lambda\in[0,1], F_1=F$ such that all $F_\lambda, \lambda\in]0,1]$ are isomonodromic and $F_0$ is a generic divide in lines. \end{definition}

It is not isomonodromic  because $F_0$ may be more degenerated than $F_1$ if it has less critical values than $F_1.$ For regular $t$ we will denote by $H_1^c(t)$ the $\C$-vector space defined by the homology of the compactification of the  fiber  $F_\lambda=t,\lambda\in[0,1]$ with coefficients in $\C$.  We prove following Theorems.

\begin{theorem}\label{simplefamily} Let $F_\lambda$ be some simple connecting family, $\lambda\in[0,1]$. Then the $\C$-vector space generated by the orbit of some oval $\delta(t)$ of the fiber $F_1=t$ contains $H_1^c(t)$.\end{theorem}

\begin{theorem}\label{structure}If $F_\lambda(x,y)$ is a simple connecting family of Hamiltonians, then for $\lambda\ne0$ there exist $\nu\leq d$ functions $\Psi_1,\cdots,\Psi_\nu$ such that each function $\Psi_k$ has logarithmic ramifications at some infinity points of the fibers $F_\lambda^{-1}(t)$ for regular values $t$ and is univalued out of the infinity points, and the \PPPF is an element of  $\C(t)[x,y,\Psi_1,\cdots,\Psi_\nu]$ .

\end{theorem}

\section{Examples}

\subsection{Perturbations of a product of 3 or 4 linear factors in general position}
The polynomial $F$ is a generic divide in lines of degree 3 or 4. We choose some oval $\delta(t)$ for regular $t$ and we denote by  ${\bf Orb}(\delta(t))$ the vector space generated by the orbit of this oval under the monodromy action.

\begin{figure}[h]
 \begin{center}
 \includegraphics[scale=0.6]{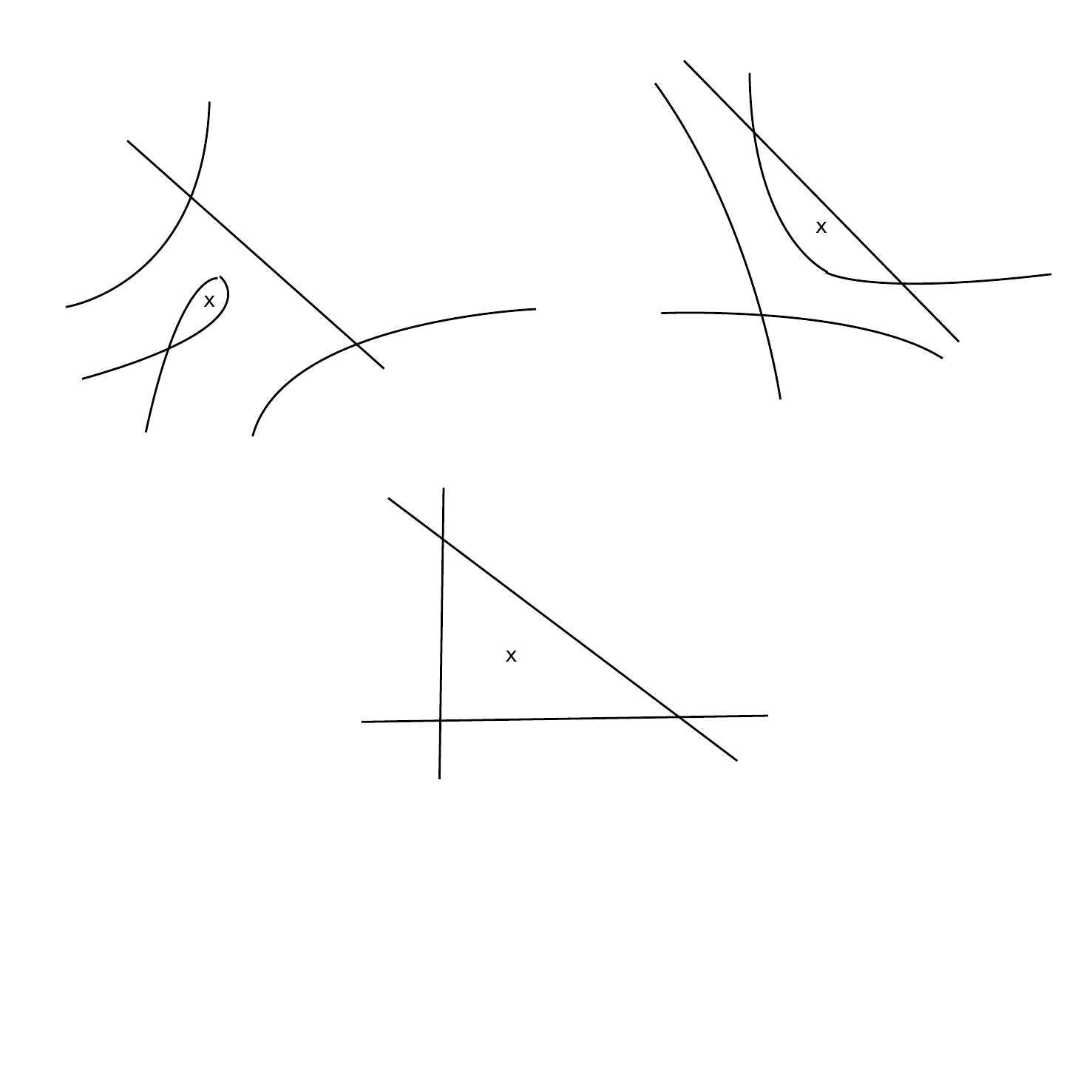} \hfill
\caption{\label{triangle}Singular points of the Hamiltonian triangle and its perturbations $(xy+\varepsilon)(x+y-1)$, for positive values of the parameter on the left and negative values on the right.}
\end{center}
\end{figure}

\medskip
If degree($F$)=3 then for regular $t,\  {\bf Orb}(\delta(t))$  is a 2-dimensional $\C$-vector space and contains the homology of the compactification of regular fibers, that is $H_1^c(t)$. It is complementary to the $\C$-vector space generated by two residual cycles at infinity. The coordinates can be chosen in such a way that $F(x,y)=xy(x+y-1)$ and $F_\varepsilon=(xy+\varepsilon)(x+y-1)$ (Figure \ref{triangle}). The 0-level and the critical points of $F$ and $F_\varepsilon$ are shown in Figure \ref{triangle}. The \PPPF is not an Abelian integral \cite{I, U}. 

The monodromy group of the Hamiltonian Triangle has only 2 generators: the monodromy around the critical value 0 and the monodromy around the center type critical value. 
Since Morse points are stable singularities, there remains three saddles and one center after perturbation if $\varepsilon$ is sufficiently small but one saddle lies on a nonzero critical level, and now the monodromy group has 3 generators. Up to reparametrization the family $F_\varepsilon$ is a simple connecting family, and $F_0$ is strictly more degenerated than $F_\varepsilon, \varepsilon\ne 0$. One can check that the codimension of ${\bf Orb}(\delta(t))$ is 1 in the perturbed case. Then the \PPPF is computed using only one multivalued function $\phi$ defined by $ d\phi=F_\varepsilon\displaystyle{d(x+y-1)\over x+y-1}=\displaystyle{(xy+\varepsilon)d(x+y-1)}$ and it is an Abelian integral.

\medskip The following  perturbation of the generic divide in 4 lines is very similar. Let
$F_0=xy\left(x+\frac12y-1\right)\left(\frac12x+y+1\right)$. The level $F=0$ contains 4 lines in general position intersecting at 6 saddle points, the homology of regular fibers is a 9 dimensional $\C$-vector space, it is generated by three vanishing cycles at center type singular values and the six cycles vanishing at 0, each one surrounding one saddle point. The only critical values are  3 critical values of center type and  0, hence the monodromy group has only 4 generators. The orbit of any oval generates a codimension 3 vector space in $H_1(t)$. Hence one needs 3 functions $\phi_1,\phi_2,\phi_3$ to perform the computation of the \PPPF (see \cite{U2} for details) and this \PPPF is generically not an Abelian integral.

\begin{figure}[h]
 \centering\hfill \includegraphics[scale=0.3]{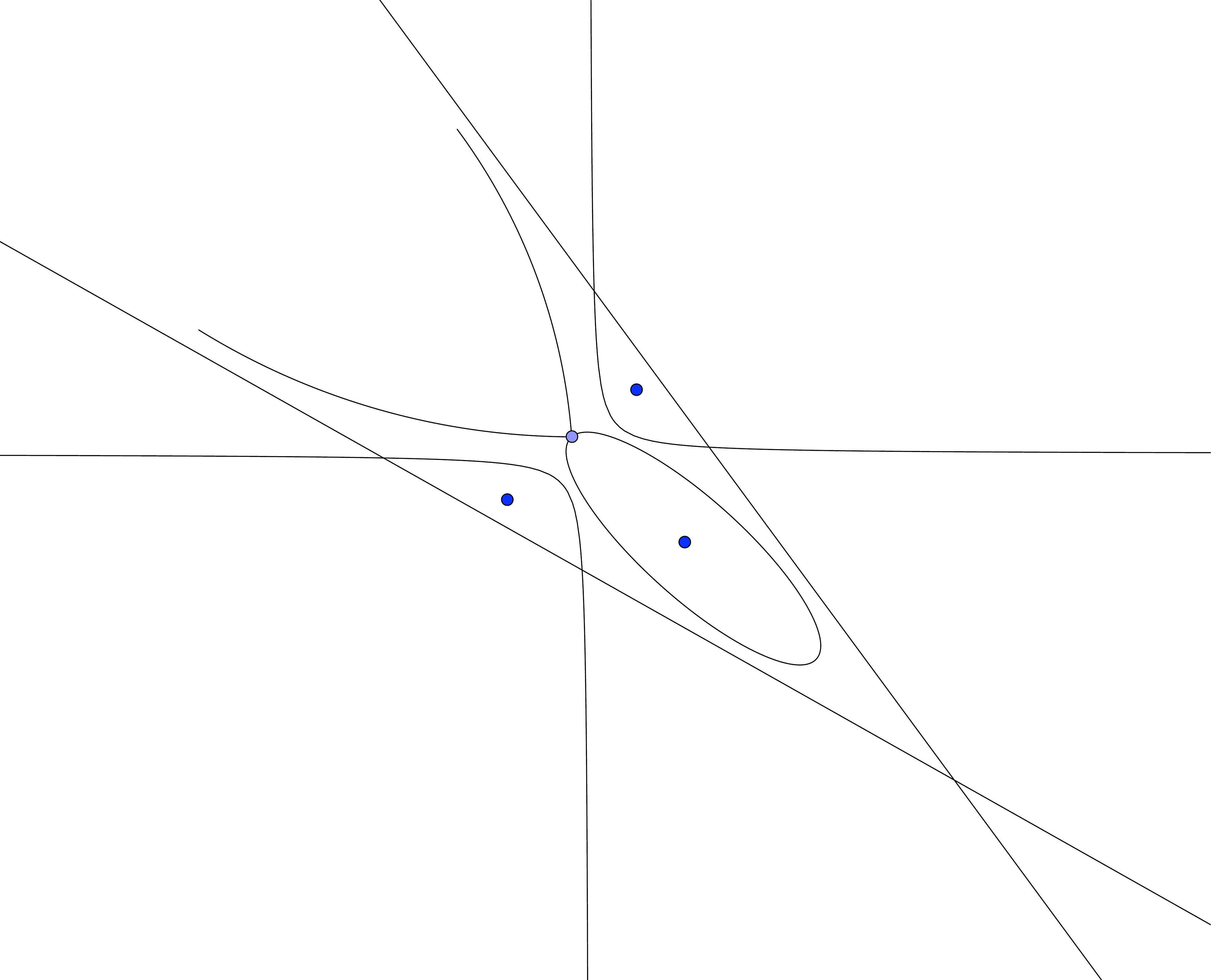} \hfill
\caption{\label{4lines}Singular curves in the phase portrait of a perturbation of 4 straight lines. This divide is included in the union of two critical levels. }
\end{figure}

 Let now 
 $$F_\varepsilon=(xy+\varepsilon)\left(x+\frac12y-1\right)\left(\frac12x+y+1\right),\varepsilon>0.$$ 
In this perturbed situation, the critical level 0 contains 5 saddle points, and the other 4 critical levels contain each one critical point, as can be seen on Figure \ref{4lines}. There are now 5 critical values. Again up to reparametrization the family $F_\varepsilon$ is a simple connecting family, and $F_0$ is strictly more degenerated than $F_\varepsilon, \varepsilon\ne 0$.
   
The orbit of any oval contains $H_1^c(t)$ and its codimension is 2. Thus one needs now only two functions to compute the Principal Poincar\'e Pontryagin Function. The \PPPF is generically not an Abelian integral, at least if the degree of the perturbative one-form is sufficiently great \cite{U2}. We can define $\phi_1,\phi_2$ as relative primitives of polynomial one-forms $F_\varepsilon\displaystyle{d\left(x+\frac12y-1\right)\over x+\frac12y-1}$ and $F_\varepsilon\displaystyle{d\left(\frac12x+y+1\right)\over \frac12x+y+1}$.
 
  We can use an additional parameter to break one more connexion:
 
 $$F_{\varepsilon,\varepsilon'}=(xy+\varepsilon)\left(\bigl( x+\frac12y-1\bigr)\bigl(\frac12x+y+1\bigr)+\varepsilon'\bigr)\right),\varepsilon'\ll\varepsilon\ll1.$$ 
 Now the codimension of ${\bf Orb}(\delta(t))$ is one and the \PPPF is an Abelian integral. We can put $F_{\varepsilon,\varepsilon'}$ into a chain of two simple connecting families, what we could call a connecting family.  Theorem \ref{simplefamily} remains true for (non necessary simple) connecting families. Note that $F_{\varepsilon,\varepsilon'}$ is not a two parameter family since we have $\varepsilon'\ll\varepsilon\ll1.$ 
 
\subsection{Codimension one case}

As usually we suppose that $F$ is a generic at infinity real Morse polynomial. 

\begin{proposition}If for regular $t$ the codimension of ${\bf Orb}(\delta(t))$ is 1 in the $\C$-vector space $H_1(t)$ then the \PPPF is an Abelian integral. \end{proposition}

\proof We denote by $r$ the dimension of  ${\bf Orb}(\delta(t))$, so that $r+1$ is the Milnor number of the fibration defined by $F$, and we denote by $\delta_1(t),\cdots,\delta_r(t)$ a basis of ${\bf Orb}(\delta(t)$ for regular $t$. We complete it to some basis $\delta_1(t),\cdots,\delta_r(t),\sigma(t)$ of $H_1(t)$, for regular $t$. The cycle $\sigma(t)$ can be chosen in such a way that it is invariant under the monodromy action thus it is a residual cycle at infinity.  
We use the monodromy representation of the \PPPF as defined in \cite{G, GI}. The same letters $\delta_k(t),\sigma(t)$ now denote loops with some base point $p(t)$ or even free loops. Let $S$ the family of loops $\delta_1,\cdots,\delta_r$. We construct a family of free loops $\hat S=\{hsh^{-1}\}, h\in \Pi_1(F^{-1}(t),p(t)), s\in S$ and finally the main geometric object $H_\delta=\displaystyle{ \hat S\over [\Pi_1(F^{-1}(t),p(t)),\hat S]}.$ This means that the elements of $H_\delta$ can be uniquely written as $\delta_1^{\alpha_1}.\cdots.\delta_r^{\alpha_r}.\sigma^0$, since we have supposed that ${\bf Orb}(\delta(t)$ does not contain $\sigma$. Thus there is a natural injection from $H_\delta$ into $H_1(F^{-1}(t),\Z)$. Moreover the group $H_\delta$ is finitely generated and from \cite{G} it has moderate growth. It only remains to use the main result of \cite{GI}. \qed 

\medskip This can be applied to the symmetric eight figure of \cite{JMP1, JMP2} or to the above mentioned  example of a perturbation of the Hamiltonian triangle. If ${\bf Orb}(\delta(t)$ contains all the homology but a 2-dimensional space generated by 2 residual cycles at infinity, again we denote by $r$ the dimension of  ${\bf Orb}(\delta(t)$ for regular $t$, now $r+2$ is the Milnor number. We complete with two cycles $\sigma_1,\sigma_2$. Then the elements of $H_\delta$ are written as, for instance,  $\delta_1^{\alpha_1}.\cdots.\delta_r^{\alpha_r}.\sigma_1^{\beta_1}.\sigma_2^{\beta_2}.\sigma_1^{\gamma_1}.\sigma_2^{\gamma_2}, \beta_1+\beta_1+\gamma_1+\gamma_1=0$ and to any such free loop there corresponds the same cycle $\alpha_1\delta_1+\cdots+{\alpha_r}\delta_r$. Clearly this map is no more injective, the \PPPF is an iterated integral on brackets and it is no more an Abelian integral, at least in the general case \cite{I, U}.

\medskip
We now assume that $F$ is a product of two irreducible polynomials $F_a,F_b$ and that the 0-level $F=0$ is a divide in the real plane and all the homology of regular fibers $F^{-1}(t)$ is seen on this divide. This is a very particular case, as we will prove. 

\begin{lemma} For any oval $\delta(t)$ the orbit ${\bf Orb}(\delta(t))$ does not contains all the homology of the fiber $F=t$.\end{lemma}

\proof The one-form $F_bdF_a$ is relatively cohomologous to $F\displaystyle{dF_a\over F_a}$. Hence its integral is not identically 0 on residual cycles around infinity points $F_a=0$. So this one-form is not algebraically relatively exact.

Nevertheless $\int_{\delta(t)}F_bdF_a=t\displaystyle{\int_{\delta(t)}{dF_a\over F_a}}\equiv 0$ since the cycle $\delta(t)$, which is vanishing at some critical value, does not turn around any point at infinity of fibers $F^{-1}(t)$.
 \qed

\begin{lemma} If the perturbation is $\omega=\alpha F_bdF_a+\beta F_adF_b$ for some reals $\alpha,\beta$, then the family of ovals is not destroyed.\end{lemma} 
 \proof First note that 
 $\omega=Fd\left(\ln(F_a^\alpha.F_b^\beta)\right)=Fdg$
 with $g=\ln(F_a^\alpha.F_b^\beta)$. The value 0 is a critical one, hence ovals are in some fiber $F=t,t\ne 0$. The only branching points of the logarithm are 0 and $\infty$. It is clear that the ovals $\delta(t)$ for regular $t$ don't turn around any branch of $F=0$, so  the function $g$ is univalued along  $\delta(t)$. The perturbed polynomial one-form is $dF+\varepsilon Fdg=d\left(F+\varepsilon Fg\right)-\varepsilon gdF.$ Since $g$ is univalued along $\delta(t),\  \int_{\delta(t)}d\left(F+\varepsilon Fg\right)\equiv0$ and of course $\int_{\delta(t)}-\varepsilon gdF\equiv0$, so that $\int_{\delta(t)}dF+\varepsilon \omega\equiv0.$ \qed
  
  \begin{figure}
\includegraphics{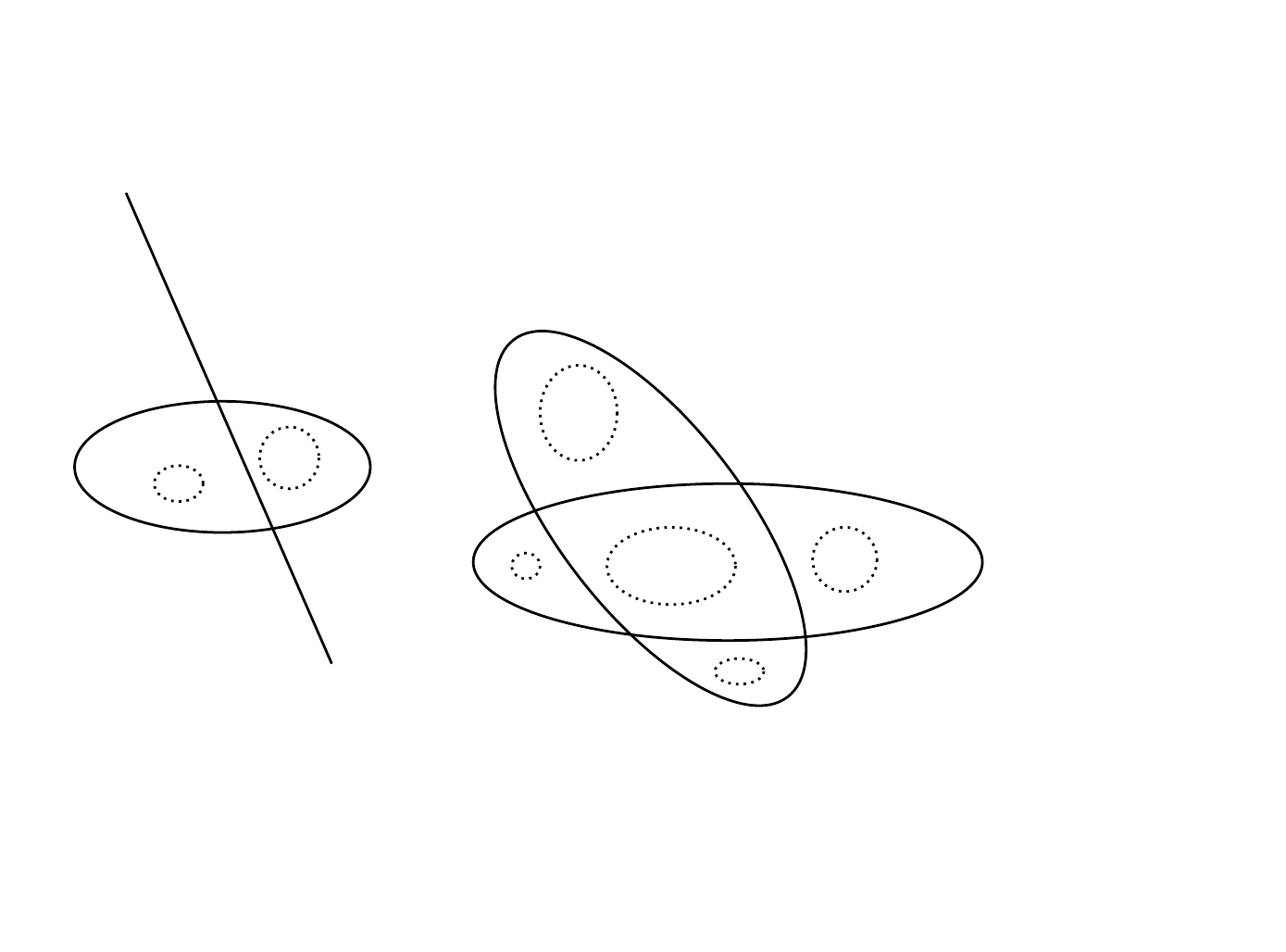}
\caption{Divide $F=0$ if $F=F_aF_b $. Dotted curves are ovals. On the left product of a line and a quadratic factor, with 2 critical points of center type and 2 double points. On the right, product of 2 quadratic factors, with 5 critical points of center type and 4 double points.\label{twofactors}}
\end{figure}
This is in fact a Darboux integrable case \cite{L}.
Following Lemma shows that this case is  very particular.

\begin{lemma} The homology is seen on the divide $F_aF_b=0$ only if the two factors are of degree at most 2 and the degree 2 factors have no real point at infinity. Then the vector space ${\bf Orb}(\delta(t)$ contains all the homology $H_1^c$ and its codimension is 1.\end{lemma}

\proof This result appears in \cite{AC}, we give a proof for completness. It is based on the computation of the Euler characteristic since the divide provides a decomposition of the disk.  On the divide we see double points and regions. We denote by $k$ the number of saddle points and $K$ the number of compact closed regions of the divide $F=0$.  We suppose that $r_\infty\leq d$ points at infinity of level curves of $F$ are real. Hence this divide has $2r_\infty$ branches going to infinity. Hence

$$1=K+2r_\infty-\left({4k\over 2}+3{2r_\infty\over 2}\right)+k+2r_\infty$$
Thus the homology is seen on the divide if

$$(d-1)^2=k+K=1+2k-r_\infty.$$
We denote by $d_a$ the degree of $F_a,$ by $d_b$ the degree of $F_b$. Using $k\leq d_ad_b$ we get

$$0\geq -r_\infty\geq d_a(d_a-2)+d_b(d_b-2).$$
The first result follows since both degrees $d_a,d_b$ are at least 1. Divides are drawn on Figure  \ref{twofactors}.

One can check that in the case of a line and an ellipse, dim$(H_1(t))=4$ and dim$({\bf Orb}(\delta(t))=3$. In the case of two ellipses, the Milnor number is 9 and dim$({\bf Orb}(\delta(t))=8$. The cycle which is not in ${\bf Orb}(\delta(t)$ has 0 as intersection number with any cycle vanishing cycle at a center critical value, thus it is a residual cycle at infinity. \qed

\section{Proofs}\label{proofs}

\subsection{Proof of Theorem \ref{simplefamily}}

In this Subsection we prove that if $F$ is in a simple connecting family $F_\lambda$, then the vector space ${\bf Orb}(\delta(t))$ contains all the homology $H_1^c(t)$ for any oval and regular values $t$. 
Therefore we prove that the monodromy group has in some sense more generators in the Milnor fibration defined by $F$ than in the Milnor fibration defined by the divide in lines $F_0$.

\begin{lemma}\label{fibration} Let $F_\lambda,\lambda\in[0,1]$ be a simple connecting family of polynomials. The map $(x,y,\lambda)\to (F_\lambda(x,y),\lambda)$ defines a fibration from $\{(x,y,\lambda)\}=\C^3$ to a subset of $\{(\lambda,t)\}=\C^2$.\end{lemma}

\proof We consider $\lambda$ as a complex parameter. From  Ehresmann Fibration Theorem \cite{W}  this mapping defines a fibration with basis the complement  of the set  where its rank is not maximal. So the basis of this fibration  is the complementary  of $\{(c,\lambda)\}$ such that $c$ is a critical value of $F_\lambda$. Hence if we denote by $\Sigma_\lambda$ the set of critical values of $F_\lambda$, the basis of this fibration is $\C^2\setminus \bigcup \{\lambda\times \Sigma_\lambda\}$. \qed

For simplicity we will project the fibers into $\{(x,y)\}=\C^2$ and denote by $F_\lambda^{-1}(t)$ any regular fiber, a Riemann surface in $\C^2$. Now again we restrict to $\lambda\in]0,1]$ where the deformation is isomonodromic. That means the following. From preceding Lemma, if $\lambda$ and $\mu$ are two near values of the parameter, there is a local connection sending $\Sigma_\lambda$ to $\Sigma_\mu$, a regular  fiber $F_\lambda^{-1}(t_\lambda)$ to some regular fiber $F_\mu^{-1}(t_\mu)$, and $H_{1,\lambda}(t_\lambda) $ onto $H_{1,\mu}(t_\mu) $. Since $F_\lambda$ and $F_\mu$ are isomonodromic they have the same number of critical values and the connection establishes a one-one corespondance  from $\Sigma_\lambda$ onto $\Sigma_\mu$. Moreover this connection commutes with the monodromy. Namely assume that some critical value $c_\lambda$ of $F_\lambda$ is sent to $C_\mu$ of $F_\mu$. Then for any cycle $\delta_\lambda(t_\lambda)$ we get the same result if we first let act the monodromy around $c_\lambda$ and then the connection, of if we first transport the cycle into $H_{1,\mu}(t_\mu)$ and the let act the monodromy around $C_\mu$. 

Now the polynomial $F_0$ is not isomonodromic with $F_\lambda$ but since all critical points are of Morse type they vary continuously with respect to $\lambda$. Hence we have the following:

\begin{lemma}The limit of critical values of $F_\lambda$ when $\lambda$ goes to 0 is one of the critical values of $F_0$.\end{lemma}

Denote as $c_j^\lambda$ the critical values of $F_\lambda$ which go to a critical value of $F_0$ of center type. Denote as $z_j^\lambda$ the critical values of $F_\lambda$ which go to 0 when $\lambda$ goes to 0.

\begin{lemma}The critical levels $c_j^\lambda$ contain only one critical point.\end{lemma}

\proof This is true for $\lambda=0$. Since it is an open property it is true for $\lambda$ near 0. By isomonodromy for $\lambda\in]0,1]$ it is true for $\lambda\in[0,1]$. \qed

The homology of the regular fibers $F_\lambda=t$ vary continuously. The vanishing cycles around critical values $c_j^\lambda$ will be denoted by $\delta_j^\lambda$. If a  cycle in the Milnor fibration defined by $F_\lambda,\lambda\in]0,1]$ vanishes at some $z_j^\lambda$ then when $\lambda\rightarrow 0$ it goes to a  vanishing cycle of the divide in lines $F_0$, more precisely a cycle vanishing at 0, which shrinks at a double point, intersection of two lines $\ell_m=0$ and $\ell_n=0$. We will denote it by $\gamma_{m,n}^\lambda$. The homology at infinity contains all residual cycles around $\ell_n=0,n=1,\cdots ,d$. Moreover we suppose that points at infinity are fixed for $\lambda\in[0,1]$. With the same notations as above, the connection of Lemma \ref{fibration} sends the homology at infinity of $F^{-1}_\lambda(t_\lambda)$ onto  the homology at infinity of $F^{-1}_\mu(t_\mu)$.

The polynomial $F_0=\ell_1.\cdots.\ell_d$ is a generic divide in lines, thus the critical level $F_0(x,y)=0$ contains $d(d-1)/2$ saddle points. Since $F_0$ is a Morse polynomial there are $(d-1)^2$ critical points, hence $K=(d-1)(d-2)/2$ critical points of center type and from genericity hypothesis all these critical points lie on distinct non zero critical levels. The monodromy operators of  the degenerated polynomial $F_0$ are the monodromy around 0 and the $K$ monodromy operators around critical values of center type, $c_j^0$. The monodromy operators of $F_\lambda, \lambda\in]0,1]$ are the monodromy around $K$ critical values $c_j^\lambda$ and  monodromy operators around each critical value $z_j^\lambda$. We will only use the monodromy generated by a loop turning once counterclockwise around all critical values $z_j^\lambda$ and only around them, and by $K$ loops turning once counterclockwise around one of the critical values $c_j^\lambda$. 

\begin{definition}The subgroup of the monodromy generated by the monodromy operators around each $c_j^\lambda$ and by a a loop turning once clockwise around all critical values $z_j^\lambda$ will be called  the sub-monodromy.\end{definition}

This sub-monodromy of $F_\lambda$  has $1+K$ generators, exactly as many generators as the monodromy of the Milnor fibration defined by $F_0$. Choose some oval $\delta(t_1)$ in the homology of one regular fiber $F_1=t_1$. It varies continuously with $\lambda$, thus it is in a family denoted by $\delta_\lambda(t_\lambda)$. Its limit when $\lambda$ goes to 0 is one of the ovals of $F_0=t_0$ where $[t_0,t_1]$ is a path in $\C$ such that $t_\lambda$ is regular for $F_\lambda,\lambda\in[0,1].$ Since everything depends continuously on $\lambda$, the orbit under the action of the sub-monodromy of $\delta(t)$ varies also continuously. Thus the dimension of the vector space generated by this orbit is constant.

\begin{notation}The dimension of ${\bf Orb}(\delta_\lambda(t))$ is denoted by $r_\lambda$.
\end{notation}

\begin{lemma}  The dimension $r_\lambda$ is at least the dimension in the degenerated case:
$r_\lambda\geq r_0$.\end{lemma}

\proof The vector space ${\bf Orb}(\delta_\lambda(t_\lambda))$ contains the vector space generated by the action of the sub-monodromy. \qed

When $\lambda$ goes to 0, the limit of ${\bf Orb}(\delta_\lambda(t_\lambda))$ is a vector space containing ${\bf Orb}(\delta_0(t_0))$ as a subspace. We know from \cite{U2} that ${\bf Orb}(\delta_0(t_0))$ contains all  homology of the compactification of regular fibers.  Moreover the homology at infinity does not depend on $\lambda$. That finishes the proof of Theorem \ref{simplefamily}.

\subsection{Proof of Theorem \ref{structure}}
If $\int_{\delta(t)}\omega$ is not identically 0 then it is an Abelian integral and it is the \PPPF and we are done. If $\int_{\delta(t)}\omega\equiv0$ we have to compute further and first we construct a convenient basis of the relative cohomology for some regular $t$.

\begin{lemma} If the orbit of some cycle $\delta(t)$ is not the whole homology, then there exist polynomial one-forms such that there integral on $\delta(t)$ is identically 0 and their integral on cycles of the complementary of ${\bf Orb} (\delta(t))$ is not 0.\end{lemma}

\proof We denote by $r$ the dimension of ${\bf Orb}(\delta(t)),\nu=(d-1)^2-r.$ The Petrov module of the integrals of polynomial one-forms on $\delta$ has dimension $r$  \cite{GaPetrov}. It contains the integrals on $\delta(t)$ of polynomial one-forms $\omega_1,\cdots,\omega_r$ free as elements of a $\C(t)$-vector space. This family can be completed with polynomial one-forms to a $\C(t)$-basis of the relative cohomology. From the dimension of the Petrov module of Abelian integrals on $\delta$, one can construct a basis of the relative cohomology of polynomial one-forms $\omega_1,\cdots,\omega_r,\psi_1,\cdots,\psi_\nu$ in such a way that the basis is  such that $\int_{\delta(t)}\psi_k\equiv 0$ if $k=1,\cdots,\nu$. \qed

\begin{remark}Since this is a basis of the relative cohomology the integrals of the one-forms $\psi_1,\cdots,\psi_\nu$ on cycles of a complementary of the orbit ${\bf Orb} (\delta(t))$ are free in the $C(t)$-vector space of Abelian integrals, hence also in the $\C[t]$-module of Abelian integrals.\end{remark}

Let us use generalized Fran\c coise's algorithm \cite{F, FP} with multivalued functions. If  $M_1(t)=-\int_{\delta(t)}\omega\equiv 0$ then there exist polynomials $ \alpha_k(F), k=1,\cdots,\nu$ such that the one-form $\omega-\sum_{k=1}^\nu\alpha_k(F)\psi_k$ has integral 0 on all cycles of $F=t$, that is this form is topologically relatively exact. From \cite{Bo} we know that there exists a polynomial $T$ in $F$ such that $T(F)\left(\omega-\sum_{k=1}^\nu\alpha_k(F)\psi_k\right)$ is algebraically relatively exact, that is there exist polynomials $Q,R$ in $x,y$ such that

$$T(F)\left(\omega-\sum_{k=1}^\nu\alpha_k(F)\psi_k\right)=Q(x,y)dF+dR(x,y).$$
This polynomial is called torsion in \cite{Bo} and it depends on reducible or non connected fibers $F=t$.
This yields

$$\omega=\sum_{k=1}^\nu\alpha_k(F)\psi_k+{dR(x,y)\over T(F)}+{Q(x,y)\over T(F)}dF.$$
Now ${dR(x,y)\over T(F)}=d\left({R(x,y)\over T(F)}\right)+T'(F){R(x,y)\over T^2(F)}dF$ where $T'$ denotes the usual derivative of the polynomial $T$ with respect to $F$.

We define a primitive of any $\psi_k,k=1,\cdots,\nu$ as follows and we will denote it by $\Psi_k$. We choose a base point in some regular fiber $F=t_0$, this allows to compute $\Psi_k$ restricted to this fiber. It is multivalued since the one-form $\psi_k$ is not algebraically relatively exact, but it is univalued along $\delta(t_0)$ and along all cycles of ${\bf Orb}(\delta(t_0)).$ Then we fix  a section transversal to our family of regular fibers. This allows to compute $\Psi_k$ on fibers $F=t$ for $t$ in our family of regular values. For any $t$, the function $\Psi_k$ is not univalued but it is univalued along all cycles of ${\bf Orb}(\delta(t)).$ Note that

$$\alpha_k(F)\psi_k=d\left(\alpha_k(F)\Psi_k\right)-\Psi_kd\left(F\alpha_k(F)\right),\ k=1,\cdots,\nu.$$

Finally there exist functions $f_1,g_1$ which are polynomials in $x,y,\Psi_1,\cdots,\Psi_\nu$ and rational in $F$ such that

$$\omega=g_1(F,x,y,\Psi_1,\cdots,\Psi_\nu)dF+df_1(F,x,y,\Psi_1,\cdots,\Psi_\nu).$$
And we can go to the next step of the algorithm and compute $M_2(t)=\int_{\delta(t)}g_1\omega$ that is a length 2 iterated integral.  This function lies in the $\C(t)$-vector space generated by Abelian integrals and integrals such as $\int_{\delta(t)}\Psi_k\psi_j$.

\begin{notation}
We denote by $I_{k,j}=\int_{\delta(t)}\Psi_k\psi_j, 1\leq k\leq \nu,1\leq j\leq\nu$. If we compute these integrals on the fibers of $F_\lambda$ we will denote them by $I_{k}^\lambda(t)$ or  $I_{k,j}^{\lambda}(t)$, and $\psi_{k}^\lambda,k=1,\cdots, \nu$ the one-forms with relative primitives $\Psi_{k}^\lambda$.
\end{notation}

Indeed since $\Psi_k\psi_j+\psi_k\Psi_j=d\left(\Psi_k\Psi_j\right)$, integrals $I_{k,j}(t)$ and $I_{j,k}\star(t)$ are opposite. So we use $I_{k,j}(t),1\leq k<j\leq r.$  Recall that from the definition of length 2 iterated integrals \cite{Ch} $I_{k,j}(t)=\int_{\delta(t)}\Psi_k\psi_j$, that is we integrate a multivalued one-form. This makes sense only on paths. So we have to suppose that there is some base point on $F=t$ and that any cycle is represented as a path. By abuse we denote again as $\delta(t)$ this path. If we change the base point the primitive $\Psi_k$ may become $\Psi_k+C_k$ for some constant $C_k$. We know from \cite{GI} that the \PPPF is base point independent.  Indeed after some change of the base point  the length-2 integral $I_{k,j}(t)$ becomes $\int_{\delta(t)}(\Psi_k+C_k)\psi_j.$ Its variation is $\int_{\delta(t)}C_k\psi_j$ which is identically 0 by construction. Thus our integrals $I_{k,j}(t)$ are base point independent. It was shown in \cite{I} that they are not Abelian. We show a little more in the following essential Lemma.

\begin{lemma}\label{free}The integrals of length 2 are free as elements of the $\C(t)$-vector space generated by Abelian integrals and the $I_{k,j}(t)$.\end{lemma}

\proof Here we consider $F=F_1$ as an element of a connecting family $F_\lambda$.
Recall that for $\lambda=0$, that is for the divide in lines, we have used one-forms $\phi_k=F{d\ell_k\over \ell_k},k=1,\cdots, d-1$. We can suppose that the complementary to ${\bf Orb}(\delta_\lambda(t))$ contains precisely the residual cycles dual to one-forms $\phi_k=F{d\ell_k\over \ell_k},k=1,\cdots,\nu$.

Furthermore the polynomials $F_\lambda$ have same points at infinity, hence same degree $d$ terms. Hence we can choose one-forms  $\psi_k,k=1,\cdots,\nu$  in such a way that for $t$ going to $\infty$ we have $\psi_k\sim\phi_k$ for $k=1,\cdots, \nu$. This yields $I_{k,j}^\lambda(t_\lambda)\sim I_{k,j}^0(t_0)$ when $t$ goes to $\infty$.

It was proved in \cite{U2} that the non Abelian integrals $I_{k,j}^0(t_0)$ are free as elements of a $\C(t)$-vector space, for $1\leq k<j\leq d-1$, hence for $1\leq k<j\leq \nu$. Now we want to know what happens if for some polynomials $\alpha_{k,j}^\lambda(F)$,

$$\sum_{k,j}\alpha_{k,j}^\lambda(t)I_{k,j}^\lambda(t)\equiv 0.$$
If these polynomials are not identically 0 they have higher degree terms. This is the dominating term of preceding combination for $t$ going to $\infty$. Since we have supposed that points at infinity of the polynomial $F_\lambda$ and one forms $\psi_k$ are independent of $\lambda$, then 
this dominating term does not depend on $\lambda$. But its limit when $\lambda$ goes to 0 is 0 because the integrals $I_{k,j,0}^\star(t)$ are free as elements of a $\C(t)$-vector space. 
Thus this dominating term is 0 and this is a contradiction.
\qed

\begin{lemma} The iterated integrals of length 3 such as
$\int_{\delta(t)}\psi_k\psi_j\psi_m$ are not base point independent.\end{lemma}

\proof Again we have to choose a base point and to integrate along paths lying in the fiber $F=t$ and starting at this base point and we denote as $\int_{\delta(t)}\psi_k\psi_j\psi_m$ an iterated integral along some path representing the cycle $\delta(t)$. If we denote by $\gamma(p)$ the piece of $\delta$ starting at the base point and going to the point $p$ of $\delta$ then $\int_{\delta(t)}\psi_k\psi_j\psi_m$ it computed as the integral along $\delta$ of the multivalued one-form which takes the value $\left(\int_{\gamma(p)}\psi_k\psi_j\right)\psi_m(p)$ at $p$ when $p$ varies along $\delta$.

The length 2 iterated integral $\int_{\gamma(p)}\psi_k\psi_j$ is computed as the integral on the path $\gamma(p)$ of the multivalued one-form $\Psi_k\psi_j$.
If we move the base point then $\Psi_k$ may become $\Psi_k+C_k$ for some constant $C_k$, the integral $\int_{\delta(t)}\psi_k\psi_j\psi_m$ becomes $\int_{\delta(t)}\psi_k\psi_j\psi_m+C_kI_{j,m}(t)$. Its variation is nonzero since this integral $I_{j,m}(t)$ is nonzero. \qed

The following result finishes the proof. 
\begin{corollary} The only base point independent one-forms are the Abelian integrals and the $I_{j,m}(t)$.\end{corollary}

It was proved in \cite{U} that generically the \PPPF of order 2, $M_2(t)$, is not an Abelian integral if the degree of the perturbative one-form is at least 5 and the Hamiltonian $F_0$ is of degree 3. 

\section{Conclusion and perspectives}A first generalization would be to use chains of simple connecting families. Namely we can use a simple connecting family connecting $F$ to a more degenerated polynomial which could be less degenerated than a generic divide in lines. The product of these paths in the manifold of Morse polynomials of degree $d$ with fixed $d$ real points at infinity is a connecting family. Proof of Theorem \ref{simplefamily} shows that the orbit of some oval in $F=t$ for generic $t$ contains all the homology $H_1^c(t)$. Theorem \ref{structure} can be generalized to such families.

Next one could relax hypothesis, for instance allow to points at infinity to move but keep generic at infinity polynomials along all the connecting family. This would allow complex points at infinity, which is natural since the monodromy works in $C^2$ and thus the fact that points at infinity are real or not is irrelevant. Therefore one has to adapt proof of Lemma \ref{free}.

We conjecture that as soon as $F$ is a Morse polynomial of degree $d$ with $d$ distinct points at infinity such that at least $d(d+1)/2$ critical levels contain only one critical point then the \PPPF is an iterated integral of length at most 2. This conjecture is based on the stratification of polynomials given by Zariski-Tarskii Theorem (or Chevalley Theorem)  \cite{B,M}. Indeed any Morse polynomial $F$ lies in some stratum. If there is a generic divide in lines in the boundary of this stratum then there exists a simple connecting family from $F$ to the generic divide in lines and we are done. Else we can use the following technic indicated to us by Maxim Kazarian who solves similar questions in \cite{K}.

We suppose that $F$ is a real polynomial of degree $d$ with $d$ real points at infinity, that all critical points are of Morse type, that at most $d/2$ critical points lie on the level $F=0$ and that other critical levels contain only one critical point. We are going to prove that such polynomials can be put in a simple connecting family.

Therefore we choose $d$ lines $\ell_1,\cdots \ell_d$ such that the algebraic curve $\ell_1\cdots\ell_d=0$ contains all $d/2$ or $(d-1)/2$ critical points of the 0-level of $F$ and both polynomials $F$ and $F_0=\ell_1\cdots\ell_d$ have the same points at infinity. Consider the family of polynomials

$$F_\lambda=\lambda F+(1-\lambda)F_0,\ F_1=F,\ F_0=\ell_1.\cdots\ell_d.$$
All polynomials $F_\lambda, \lambda\in[0,1]$ have the same points at infinity. Moreover all $d/2$ or $(d-1)/2$ critical points of the critical level $F=0$ are critical points of $F_\lambda,\lambda\in[0,1]$.

If all critical points of $F_\lambda$ are of Morse type and all non zero critical levels of $F_\lambda,\lambda\in]0,1]$ contain only one critical point, we have constructed a simple connecting family and the result is proved. Else it means that for some isolated values of $\lambda$ either one critical point is not Morse or one nonzero critical level contains more than one critical point.  By Zariski-Tarski Theorem each of these two conditions define an algebraic set in $\C_{d-1}(x,y]$, vector space of polynomials of degree at most  $d$ with fixed points at infinity.  Thus it can be avoided by following a path in $\C_{d-1}[x,y]$ instead of $\R_{d-1}[x,y]$. Again we have constructed a simple connecting family and the result is proved.

\begin{center}Ê\textsc{Acknowledgments}\end{center}This work was partially supported by
the Fondecyt Project 11080250 of the Chilian Government.

The authors thank greatly Maxim Kazarian and Pavao Marde\v si\'c for fruitful discussions, and J.-P. Rolin for his help.

\end{document}